\documentclass[letterpaper, 10 pt, conference]{ieeeconf}

\pdfoutput=1

\IEEEoverridecommandlockouts
\usepackage{amsmath} 
\usepackage{amssymb}  
\usepackage[utf8]{inputenc}
\usepackage{algorithm2e}
\usepackage{color}
\usepackage{tikz}
\usetikzlibrary{shapes,arrows.meta}

\usepackage{graphicx}
\usepackage[font=footnotesize,labelfont=bf]{caption}

\DeclareMathOperator*{\argmin}{arg\,min}

\newtheorem{remark}{Remark}
\newtheorem{prop}{Proposition}

\newtheorem{assumption}{Assumption}

\newcommand{\norm}[1]{\left\lVert#1\right\rVert}

\def\oprocend{\hfill$\square$}
\def\oprocendbis{\hfill$\blacksquare$}

\title{\LARGE \bf
A Partition-Based Implementation of the Relaxed ADMM for Distributed Convex Optimization over Lossy Networks
}
\author{N. Bastianello$^\dagger$, M. Todescato$^\ddagger$, R. Carli$^\dagger$, L. Schenato$^\dagger$
\thanks{$^\dagger$ Department of Information Engineering (DEI), University of Padova, Italy. {\tt\small nicola.bastianello.3@studenti.unipd.it, [carlirug|schenato]@dei.unipd.it}.}
\thanks{$^\ddagger$ Bosch Center for Artificial Intelligence. Renningen, Germany. {\tt\small mrc.todescato@gmail.com}. Part of the work was carried out during the author's postdoctoral fellowship at DEI.}
}

\begin{document}

\maketitle
\thispagestyle{empty}
\pagestyle{empty}

\begin{abstract}
%
In this paper we propose a distributed implementation
of the relaxed Alternating Direction Method of Multipliers algorithm
(R-ADMM) for optimization of a separable convex cost
function, whose terms are stored by a set of interacting agents,
one for each agent. Specifically  the local cost stored by each node is in
general a function of both the state of the node and the states of its
neighbors, a framework that we refer to as `partition-based' optimization.
This framework presents a great flexibility and can be adapted to a large
number of different applications. We
show that the partition-based R-ADMM algorithm we introduce
is linked to the relaxed Peaceman-Rachford Splitting (R-PRS)
operator which, historically, has been introduced in the literature to
find the zeros of sum of functions.
Interestingly, making use of non expansive operator theory, the proposed
algorithm is shown to be provably robust against random packet losses that
might occur in the communication between
neighboring nodes. Finally, the effectiveness of the proposed algorithm is
confirmed by  a set of compelling numerical simulations run over random
geometric graphs subject to i.i.d. random packet losses.
%
\end{abstract}


\section{Introduction}\label{sec:intro}
Because of the advent of the Internet-of-Things (IoT), we are witnessing the proliferation of large-scale systems in which a multitude of locally networked peers are able to asynchronously and unreliably stream information with neighboring peers. In such networks, many system-wide applications, e.g., in the machine learning field \cite{Slavakis:14}, as well as operational requirements, e.g., state estimation and control, can be cast as optimization problems aiming at globally optimal configurations across the entire network. Yet, many of these applications, because of the locally connected nature of such systems, are characterized by local decision functions which are influenced only by local information flows among neighboring peers. Examples owning to the above framework can be found in applications such as state estimation and power flow control in Smart Electric Grids \cite{todescato2017generalized} as well as cooperative localization in Wireless Networks \cite{carli2013distributed}, just to mention a few. Such class of problems, referred to as \emph{partition-based} optimization problems, can be formally described as
\begin{equation}\label{eq:partion-based-problem}
\min_{x_i, i\in\mathcal{V}} \sum_{i=1}^N f_i\left(x_i,\{x_j\}_{j\in\mathcal{N}_i}\right)
\end{equation}
where each peer $i$, in the set $\mathcal{V}=\{1,\ldots, N\}$ of all possible peers, is responsible for a local decision function $f_i$ only affected by its local piece of information $x_i$ as well as by $x_j$ of its neighboring peers $j\in\mathcal{N}_i$.\\
Regarding partioned-based optimization, many different approaches have been analyzed in the recent literature. In \cite{carli2013distributed} a dual decomposition approach is proposed. Gradient-based schemes are considered in \cite{todescato2015robust,todescato2017generalized}. In \cite{todescato2015robust}, the authors present a Block-Jacobi iteration suitable for quadratic programming while in the more recent \cite{todescato2017generalized} a modified generalized scheme is presented for generic convex optimization.
Other solutions involve the well-know Alternating Direction Method of Multipliers (ADMM) \cite{glowinski1975approximation,gabay1976dual} which has been shown to be particularly suited for parallel and distributed computations. We refer the reader to \cite{boyd2011distributed,fukushima1992application,eckstein1994some,eckstein1992douglas,chen1994proximal,EG:AT:ES:MJ:2015}
for an overview of possible applications and convergence results. In particular, to the best of the authors knowledge, a partioned-based ADMM scheme has been first introduced in \cite{erseghe2012distributed} to solve for cooperative localization in WSN, while in \cite{mota2015distributed} partition-based ADMM is applied to MPC.\\
However, when dealing with asynchronous and, in particular, faulty/unreliable communications, i.e., subject to delays and packet drops, while first-order \cite{todescato2015robust,todescato2017generalized} and second-order \cite{carli2015distributed,carli2015analysis} gradient-based schemes have been proved to be robust, the same results do not apply to ADMM schemes. More specifically, works devoted to the study of asynchronous ADMM  implementations can be found. Examples are \cite{bianchi2016coordinate}, where convergence of the ADMM is shown when only a subset of coordinates is randomly updated at every time instant, and the recent \cite{peng2016arock}, where a framework for asynchronous operations is proposed. Conversely, literature on robustness analysis of ADMM in the presence of faulty communication is still scarce and usually confined to specific setups such as \cite{zhang2014asynchronous,chang2016asynchronous} where only bounded delays and a particular master-slave communication architecture are considered. To the authors knowledge, first steps toward more general results have appeared only recently in \cite{bastianello2018distributed} where a robust generalized ADMM scheme for consensus optimization is presented.\\
In this paper we take over from \cite{bastianello2018distributed}. We leverage nonexpansive operator theory where the underlying idea is to reformulate the original optimization problem into an equivalent form whose solutions corresponds to the fixed points of a suitable operator. One particular class of such operators is represented by \emph{splitting methods} in which the burden of solving for the fixed points is alleviated by breaking down the computations into several steps. Well-known examples are the Peaceman-Rachford Splitting (PRS) \cite{peaceman1955numerical} with its generalized versions as well as the Douglas-Rachford Splitting (DRS) \cite{douglas1956numerical,lions1979splitting}. Moreover we refer to \cite{davis2016convergence,hannah2016unbounded} for more details and possible applications to asynchronous setups.\\
We start our analysis from the fact that ADMM can be shown to be equivalent to the DRS applied to the Lagrange dual of the original problem \cite{eckstein2012augmented}. Then, the contribution of the paper are twofold. First, we present a reformulated version of ADMM suitable for partition-based optimization. Second, by resorting to results on stochastic operator theory, we formally prove robustness of the proposed algorithm to faulty communications.\\
The remainder of the paper is organized as follows. Section~\ref{sec:introADMM} reviews the classical ADMM algorithm and its generalized version, referred to as Relaxed-ADMM. Section~\ref{sec:distributed_consensus} introduces the partition-based framework for consensus optimization and derives the relaxed ADMM applied to this problem. Section \ref{sec:robustADMM} describes the proposed robust implementation. Section~\ref{sec:simulation} collects some numerical simulations. Finally, Section~\ref{sec:conclusions} draws some concluding remarks. Due to space constraints, all the technical proofs can be found in the Appendices.
\section{The Relaxed-ADMM algorithm}\label{sec:introADMM}

Consider the following optimization problem
\begin{align}\label{eq:primal-problem}
\begin{split}
	&\min_{x\in\mathcal{X},y\in\mathcal{Y}} \{f(x)+g(y)\}\\
	&\text{s.t.}\ Ax+By=c
\end{split}
\end{align}
with $\mathcal{X}$ and $\mathcal{Y}$ Hilbert spaces, $f:\mathcal{X}\rightarrow\mathbb{R}\cup\{+\infty\}$ and $g:\mathcal{Y}\rightarrow\mathbb{R}\cup\{+\infty\}$ closed, proper and convex functions\footnote{A function $f:\mathcal{X}\rightarrow\mathbb{R}\cup\{+\infty\}$ is said to be \textit{closed} if $\forall a \in\mathbb{R}$ the set $\{x\in\operatorname{dom}(f)\ |\ f(x)\leq a\}$ is closed. Moreover, $f$ is said to be \textit{proper} if it does not attain $-\infty$ \cite{boyd2011distributed}.}. In the following we assume that the above problem has solution.

Let us define the \textit{augmented Lagrangian} for the problem \eqref{eq:primal-problem} as
\begin{align}\label{eq:augmented-lagr}
	\mathcal{L}_{\rho}(x,y;w)=&f(x)+g(y)-w^\top\left(Ax+By-c\right)\nonumber\\
	&+\frac{\rho}{2}\|Ax+By-c\|^2
\end{align}
where $\rho>0$ and $w$ is the vector of Lagrange multipliers. 

The Relaxed-ADMM (R-ADMM) algorithm (see \cite{davis2016convergence}) consists in the alternating of the following three steps
\begin{align}
\begin{split}
	y(k+1)&=\argmin_y\{\mathcal{L}_\rho(x(k),y;w(k))\\&+\rho(2\alpha-1)\langle By,(Ax(k)+By(k)-c)\rangle\}\label{eq:r-admm-1}
\end{split}\\
\begin{split}
	w(k+1)&=w(k)-\rho(Ax(k)+By(k+1)-c)\\&-\rho(2\alpha-1)(Ax(k)+By(k)-c)\label{eq:r-admm-2}
\end{split}\\
	x(k+1)&=\argmin_x\mathcal{L}_\rho(x,y(k+1);w(k+1)).\label{eq:r-admm-3}
\end{align}

The R-ADMM algorithm can be derived applying the relaxed \textit{Peaceman-Rachford splitting operator} to the Lagrange dual of problem \eqref{eq:primal-problem} \cite{boyd2011distributed,davis2016convergence}. It can be shown that, under the assumptions made on functions $f$ and $g$, the convergence of the R-ADMM algorithm is guaranteed if 
\begin{align*}
0< \alpha < 1, \qquad \rho>0.
\end{align*}
We conclude this section by observing that setting $\alpha=1/2$ one can retrieve the classical ADMM algorithm widely analyzed in \cite{boyd2011distributed}.

\section{Distributed Partition-Based Convex Optimization}\label{sec:distributed_consensus}

\subsection{Problem Formulation}\label{subsec:distributed_problem}

We start by formulating the problem we aim at solving.

Let $\mathcal{G}=(\mathcal{V},\mathcal{E})$ be an undirected graph, where $\mathcal{V}$ denotes the set of $N$ vertices, labeled $1$ through $N$, and $\mathcal{E}$ the set of edges. For $i \in \mathcal{V}$, by $\mathcal{N}_i$ we denote the set of neighbors of node $i$ in $\mathcal{G}$, namely,
$$
\mathcal{N}_i =\left\{j \in \mathcal{V} \,:\, (i,j) \in \mathcal{E} \right\}.
$$
The state of each node is characterized by the local variable $x_i\in\mathbb{R}^n$, and we are interested in solving the following optimization problem
\begin{align}\label{eq:opt_problem}
\begin{split}
	&\min_{x_i,\ i\in\mathcal{V}}\sum_{i=1}^Nf_i\left(x_i,\{x_j\}_{j\in\mathcal{N}_i}\right)
\end{split}
\end{align}
where $f_i: \mathbb{R}^{n|\mathcal{N}_i|} \rightarrow\mathbb{R}\cup\{+\infty\}$ are closed, proper and convex functions and where $f_i$ is known only to node $i$\footnote{Note that in general the local variables might have different dimensions, but here for simplicity they are assumed to be all vectors in $\mathbb{R}^n$.}. Note that writing $f_i\left(x_i,\{x_j\}_{j\in\mathcal{N}_i}\right)$ denotes that $f_i$ is a function of node $i$'s state and the states of its neighbors only. We refer to this framework as to \emph{partition-based convex optimization}.

In the following we denote by $\mathbf{x}^*\in\mathbb{R}^{nN}$ the optimal solution of \eqref{eq:opt_problem} with components $x_i^*\in\mathbb{R}^n$.

In this work we seek for iterative and distributed algorithms solving problem in \eqref{eq:opt_problem} where, by distributed, we mean that each node can exchange information only with its neighbors. 
To this goal, we provide an alternative formulation of \eqref{eq:opt_problem} to which we can apply the relaxed ADMM reviewed in the previous Section. 

As first step, we assume that each node $i$ stores local copies of the variables its cost function depends on, denoted by
\begin{equation*}
	x_i^{(i)}\ \text{and}\ x_j^{(i)}\ \forall j\in\mathcal{N}_i;
\end{equation*}
precisely, $x_j^{(i)}$ is the local copy of the variable $x_j$ stored in memory by node $i$, while $x_i^{(i)}=x_i$. Observe that, problem in \eqref{eq:opt_problem} can be equivalently formulated as
\begin{align}\label{eq:distributed-primal}
\begin{split}
	&\min_{x_i^{(i)},\forall i}\sum_{i=1}^Nf_i\left(x_i^{(i)},\{x_j^{(i)}\}_{j\in\mathcal{N}_i}\right)\\
	&\text{s.t.}\ \ x_i^{(i)}=x_i^{(j)},\\
	&\qquad x_j^{(i)}=x_j^{(j)},\ \forall (i,j)\in\mathcal{E}
\end{split}
\end{align}
where the constraints impose that the consensus be reached, that is, all local copies of each vector $x_i$ be equal.

Now, for each edge $(i,j)\in\mathcal{E}$, we introduce the \textit{bridge variables} $y_i^{(i,j)}$ and $y_j^{(i,j)}$, and $y_i^{(j,i)}$ and $y_j^{(j,i)}$. 
Notice that the constraints in \eqref{eq:distributed-primal} can be rewritten as
\begin{align}\label{eq:EquivalentConstraints}
\begin{split}
	& x_i^{(i)}=y_i^{(i,j)}, \quad x_j^{(i)}=y_j^{(i,j)} \\
	& x_i^{(j)}=y_i^{(j,i)}, \quad x_j^{(j)}=y_j^{(j,i)} \\
	& y_i^{(i,j)}=y_i^{(j,i)}, \quad y_j^{(i,j)}=y_j^{(j,i)}
\end{split}\ \ \ \forall (i,j)\in\mathcal{E}.
\end{align}
We define now the vectors
\begin{equation*}
	\mathbf{x}^{(i)} = \begin{bmatrix}
		x_i^{(i)} \\ \{x_j^{(i)}\}_{j\in\mathcal{N}_i}
	\end{bmatrix}
	\quad \text{and}\quad
	\mathbf{y}^{(i)} = \begin{bmatrix}
		\{y_i^{(i,j)}\}_{j\in\mathcal{N}_i} \\ \{y_j^{(i,j)}\}_{j\in\mathcal{N}_i}
	\end{bmatrix}
\end{equation*}
where $\mathbf{x}^{(i)}\in\mathbb{R}^{n(|\mathcal{N}_i|+1)}$ and $\mathbf{y}^{(i)}\in\mathbb{R}^{2n|\mathcal{N}_i|}$, and the overall vectors $\mathbf{x}=[\mathbf{x}^{(1)\top},\ldots,\mathbf{x}^{(N)\top}]^\top$ and $\mathbf{y}=[\mathbf{y}^{(1)\top},\ldots,\mathbf{y}^{(N)\top}]^\top$. Let $f(\mathbf{x})=\sum_{i=1}^N f_i(\mathbf{x}^{(i)})$, then problem \eqref{eq:distributed-primal} with constraints \eqref{eq:EquivalentConstraints}, 
can be compactly rewritten as
\begin{align*}
	& \min_{\mathbf{x}} f(\mathbf{x})\\
	& \text{s.t.}\ \ A\mathbf{x}+\mathbf{y}=0\\
	&\qquad  \mathbf{y}=P\mathbf{y}
\end{align*}
for a suitable $A$ matrix and with $P$ being a permutation matrix that swaps $y_i^{(i,j)}$ with $y_i^{(j,i)}$. Making use of the indicator function $\iota_{(I-P)}(\mathbf{y})$ which is equal to 0 if $(I-P)\mathbf{y}=0$, and $+\infty$ otherwise, we can finally rewrite problem \eqref{eq:distributed-primal} as
\begin{align}\label{eq:primal-indicator-f}
\begin{split}
	& \min_{\mathbf{x},\mathbf{y}}\left\{f(\mathbf{x})+\iota_{(I-P)}(\mathbf{y})\right\}\\
	& \text{s.t.}\ \ A\mathbf{x}+\mathbf{y}=0.
\end{split}
\end{align}
Clearly problem \eqref{eq:primal-indicator-f} conforms to the formulation of problem \eqref{eq:primal-problem}, therefore we can apply the R-ADMM to solve it, which is the focus of the following Section.

\subsection{R-ADMM for Partition-Based Optimization}\label{subsec:distributed_R-ADMM}
In this Section we propose an implementation of the R-ADMM for the partition-based convex optimization problems introduced in Section \ref{subsec:distributed_problem}.\\

It is possible to show that
the direct application of the R-ADMM algorithm to \eqref{eq:primal-indicator-f} is amenable of distributed implementation. More precisely,
let $w_i^{(i,j)}$ be the Lagrange multiplier associated to the constraint $x_i^{(i)}=y_i^{(i,j)}$ and let
$$
	\mathbf{w}^{(i)} = \begin{bmatrix}
		\{w_i^{(i,j)}\}_{j\in\mathcal{N}_i} \\ \{w_j^{(i,j)}\}_{j\in\mathcal{N}_i}
	\end{bmatrix}.
$$
Assume that node $i$ stores in memory $\mathbf{x}^{(i)}$, $\mathbf{y}^{(i)}$ and $\mathbf{w}^{(i)}$. Then one can see that these variables can be updated by node $i$ (according to the updating equations of R-ADMM) only receiving the information 
$\mathbf{x}^{(j)}$, $\mathbf{y}^{(j)}$ and $\mathbf{w}^{(j)}$, for all $j \in \mathcal{N}_i$; namely, only communications among neighboring nodes are required. We do not report the explicit equations, because they are quite unwieldy.

However, leveraging the particular structure of problem \eqref{eq:primal-indicator-f} and considering the fact that we are interested only into the trajectories of the variables $\mathbf{x}^{(i)}$, it is possible to provide a simpler implementation of the R-ADMM algorithm, in terms of memory and communication requirements; in particular, this lighter implementation, beside the variables $\mathbf{x}^{(i)}$ involve only the auxiliary variables $z_i^{(i,j)}$, $i \in \mathcal{V}$, $j \in \mathcal{N}_i$.

We have the following Proposition.


\smallskip

\begin{prop}\label{prop:distributed-r-admm}
The trajectories generated by the variables $\mathbf{x}^{(i)}$, $i \in V$, obtained by applying the R-ADMM algorithm in \eqref{eq:r-admm-1}, \eqref{eq:r-admm-2}, \eqref{eq:r-admm-3}, to the problem in \eqref{eq:primal-indicator-f}, starting from a given initial condition $\mathbf{x}^{(i)}(0)$, $\mathbf{w}^{(i)}(0)$, $\mathbf{y}^{(i)}(0)$ are identical to the trajectories generated by iterating the following equations
\begin{align}\label{eq:x-update-distributed}
	\mathbf{x}^{(i)}(k) &= \argmin_{x_i^{(i)},\{x_j^{(i)}\}_{j\in\mathcal{N}_i}} \Bigg\{ f_i\left(x_i^{(i)},\{x_j^{(i)}\}_{j\in\mathcal{N}_i}\right)+ \nonumber\\
	& -\left(\sum_{j\in\mathcal{N}_i}z_i^{(j,i)}(k)\right)^\top x_i^{(i)}-\sum_{j\in\mathcal{N}_i}z_j^{(j,i)}(k)^\top x_j^{(i)}+ \nonumber \\
	& \left. +\frac{\rho}{2}|\mathcal{N}_i|\|x_i^{(i)}\|^2+\frac{\rho}{2}\sum_{j\in\mathcal{N}_i}\|x_j^{(i)}\|^2 \right\}
\end{align}
for all $i \in V$, and 
\begin{align}\label{eq:z-update-distributed}
\begin{split}
	& z_i^{(i,j)}(k+1) = (1-\alpha)z_i^{(i,j)}(k)-\alpha z_i^{(j,i)}(k)+2\alpha\rho x_i^{(i)}(k) \\
	& z_j^{(i,j)}(k+1) = (1-\alpha)z_j^{(i,j)}(k)-\alpha z_j^{(j,i)}(k)+2\alpha\rho x_j^{(i)}(k) \\
\end{split}
\end{align}
for all $j\in\mathcal{N}_i$, where the auxiliary variables are initialized as $z_i^{(i,j)}(0)= w_i^{(i,j)}(0)  + \rho y_i^{(i,j)}(0)$.

\oprocend
\end{prop}

\smallskip

The proof of the previous Proposition can be found in Appendix~\ref{app:distributed-admm}.

Proposition \ref{prop:distributed-r-admm} suggests a straightforward distributed implementation of the R-ADMM in which a node $i$ locally stores and updates the variables $x_i^{(i)}$ and $x_j^{(i)}$, $z_i^{(i,j)}$, $z_j^{(i,j)}$ for all $j\in\mathcal{N}_i$. Within this implementation the node requires the auxiliary variables $z_i^{(j,i)}$ and $z_j^{(j,i)}$ to be sent by each of its neighbors in order to update the local $x$ variables and hence the auxiliary variables.\\
An equivalent implementation can be obtained if node $i$ stores and updates the $z_i^{(j,i)}$ and $z_j^{(j,i)}$ variables instead of the $z_i^{(i,j)}$ and $z_j^{(i,j)}$ variables. This implementation is formally described in Algorithm \ref{alg:smart-distributed-r-admm}.

\smallskip

\begin{algorithm}[ht!]
	\SetKwInOut{Input}{Input}
	\Input{Set the termination condition $K>0$. For each node $i$, initialize $\mathbf{x}^{(i)}(0)$ and $\{z_i^{(j,i)}(0),z_j^{(j,i)}(0)\}_{j\in\mathcal{N}_i}$.}
	$k\leftarrow0$\;
	\While{$k<K$ each agent $i$}{
		compute $\mathbf{x}^{(i)}(k)$ according to \eqref{eq:x-update-distributed}\;
		for all $j \in \mathcal{N}_i$, compute the temporary variables 
		\begin{align}\label{eq:q}
		\begin{split}
			& q_i^{(i \to j)} = -z_i^{(j,i)}(k)+2 \rho x_i^{(i)}(k) \\
			& q_j^{(i \to j)} = -z_j^{(j,i)}(k)+2 \rho x_j^{(i)}(k)
		\end{split};
		\end{align}\\
		for all $j \in \mathcal{N}_i$, transmit $\{q_i^{(i \to j)},q_j^{(i \to j)}\}$ to node $j$\;
		gather $\{q_j^{(j \to i)},q_i^{(j \to i)}\}$ from each neighbor $j$\;
		update the auxiliary variables as
		\begin{align}\label{eq:Update_z_ji}
		\begin{split}
			& z_i^{(j,i)}(k+1) = (1-\alpha)z_i^{(j,i)}(k)+ \alpha q_i^{(j \to i)} \\
			& z_j^{(j,i)}(k+1) = (1-\alpha)z_j^{(j,i)}(k)+ \alpha q_j^{(j \to i)}
		\end{split};
		\end{align}\\
		$k\leftarrow k+1$\;
	}
\caption{Modified partition-based R-ADMM.}
\label{alg:smart-distributed-r-admm}
\end{algorithm}

\smallskip

Observe, that, at the beginning of each iteration node $i$ updates $\mathbf{x}^{(i)}$ based only on local information. Then it computes the temporary variables $q_i^{(i \to j)}$ and $q_j^{(i \to j)}$ which are sent to neighbor $j$. At the same time, it receives the quantities 
 $q_j^{(j \to i)}$ and $q_i^{(j \to i)}$ from neighbor $j$ and it uses these information to update $z_i^{(j,i)}$ and $z_j^{(j,i)}$ as in \eqref{eq:Update_z_ji}. 

The following Proposition characterizes the convergence properties  of Algorithm~\ref{alg:smart-distributed-r-admm}, which follows from those of the R-PRS. The proof is available in Appendix~\ref{app:convergence}.

\smallskip

\begin{prop}\label{prop:convergence}
For Algorithm \ref{alg:smart-distributed-r-admm} let $(\alpha, \rho)$ be such that $0<\alpha <1$ and $\rho >0$. Then, for any initial conditions, the trajectories $k \to x_i^{(i)}(k)$ and $k \to x_j^{(i)}(k)$, $i \in \mathcal{V}$, $j\in\mathcal{N}_i$, generated by the Algorithm~\ref{alg:smart-distributed-r-admm}, converge to the optimal solution of \eqref{eq:opt_problem}, i.e.,
\begin{align*}
\begin{split}
	\lim_{k \to \infty} x_i^{(i)}(k) = x_i^* \\
	\lim_{k \to \infty} x_j^{(i)}(k) = x_j^* 
\end{split}, \qquad \forall i \in \mathcal{V}, j\in\mathcal{N}_i.
\end{align*}
\oprocend
\end{prop}

\begin{remark}
The implementation of the R-ADMM presented in Algorithm \ref{alg:smart-distributed-r-admm} requires that each node stores and updates locally $3|\mathcal{N}_i|+1$ variables. Moreover, a node has to transmit only two variables to each of its neighbors at each time instant.
\end{remark}

\section{Partition-based R-ADMM over lossy networks}\label{sec:robustADMM}
The partition-based algorithm described in the previous Section is proved to converge under the implicit assumption that the communication channels are reliable and, therefore, no packet loss occurs.\\
The aim of this Section is to prove the convergence of the partition-based R-ADMM in case the communications are unreliable and some transmissions between nodes might fail. Precisely, we make the following assumption.
\smallskip

\begin{assumption}\label{ass:lossy}
During any iteration of Algorithm \ref{alg:smart-distributed-r-admm}, the communication from node $i$ to node $j$ can be lost with some probability $p$.\oprocend
\end{assumption}

\smallskip

To formally describe the communication failures, we associate to each transmission a random variable which is equal to $1$ if the packet is lost, $0$ otherwise. Specifically we introduce the family of binary random variables $L^{(i \to j)}(k)$, $k=0,1,2,\ldots$, $i \in \mathcal{V}$, $j \in \mathcal{N}_i$, which are independent for $i,\ j$ and $k$ that vary, such that
$$
\mathbb{P}\left[L^{(i \to j)}=1\right]=p, \qquad \mathbb{P}\left[L^{(i \to j)}=0\right]=1-p.
$$
%
Accounting for the potential packet losses that we have introduced above, Algorithm \ref{alg:smart-distributed-r-admm} is modified as illustrated in \ref{alg:robust-smart-distributed-r-admm}.

\smallskip

\begin{algorithm}[ht!]
	\SetKwInOut{Input}{Input}
	\Input{Set the termination condition $K>0$. For each node $i$, initialize $\mathbf{x}^{(i)}(0)$ and $\{z_i^{(j,i)}(0),z_j^{(j,i)}(0)\}_{j\in\mathcal{N}_i}$.}
	$k\leftarrow0$\;
	\While{$k<K$ each agent $i$}{
		compute $\mathbf{x}^{(i)}(k)$ according to \eqref{eq:x-update-distributed}\;
		compute, for $j \in \mathcal{N}_i$, the temporary variables 
		\begin{align*}
		\begin{split}
			& q_i^{(i \to j)} = -z_i^{(j,i)}(k)+2 \rho x_i^{(i)}(k) \\
			& q_j^{(i \to j)} = -z_j^{(j,i)}(k)+2 \rho x_j^{(i)}(k)
		\end{split};
		\end{align*}\\
		transmit, for $j \in \mathcal{N}_i$, $\{q_i^{(i \to j)},q_j^{(i \to j)}\}$ to node $j$\;
		\For{$j\in\mathcal{N}_i$, if $\{q_j^{(j \to i)},q_i^{(j \to i)}\}$ is received}{
			update the auxiliary variables as 
			\begin{align*}
			\begin{split}
				& z_i^{(j,i)}(k+1) = (1-\alpha)z_i^{(j,i)}(k)+ \alpha q_i^{(j \to i)} \\
				& z_j^{(j,i)}(k+1) = (1-\alpha)z_j^{(j,i)}(k)+ \alpha q_j^{(j \to i)}
			\end{split};
			\end{align*}
		}
		$k\leftarrow k+1$\;
	}
\caption{Robust partition-based R-ADMM.}
\label{alg:robust-smart-distributed-r-admm}
\end{algorithm}

\smallskip

\noindent In this potentially lossy scenario, during the $k$-th iteration, node $i$ updates the local variables $x_i^{(i)}$ and $\{x_j^{(i)}\}_{j\in\mathcal{N}_i}$, according to \eqref{eq:x-update-distributed}. Then, it computes the temporary variables $\{q_i^{(i \to j)},q_j^{(i \to j)}\}$ for each of its neighbors $j\in\mathcal{N}_i$ as in \eqref{eq:q} and transmits them. If neighbor $j$ receives the packet, that is, if $L^{(i \to j)}(k)=0$, then it updates the auxiliary variables $z_j^{(i,j)}$ and $z_i^{(i,j)}$ using the received values $\{q_i^{(i \to j)},q_j^{(i \to j)}\}$ according to \eqref{eq:Update_z_ji}, otherwise it leaves them unchanged.\\
The updates for the auxiliary variables can be described in a compact way as follows
\begin{align*}
	z_i^{(j,i)} & (k+1) = L^{(j \to i)}(k)z_i^{(j,i)}(k) + \\
		& + \left(1-L^{(j \to i)}(k)\right) \, \left( (1-\alpha)z_i^{(j,i)}(k)+\alpha q_i^{(j\to i)}\right) \\
	z_j^{(j,i)} & (k+1) = L^{(j \to i)}(k)z_j^{(j,i)}(k) + \\
		& + \left(1-L^{(j \to i)}(k)\right) \, \left( (1-\alpha)z_j^{(j,i)}(k)+\alpha q_j^{(j\to i)}\right).
\end{align*}

The following Proposition characterizes the convergence properties of Algorithm \ref{alg:robust-smart-distributed-r-admm}, running in the probabilistic lossy scenario described in Assumption \ref{ass:lossy}. 

\smallskip

\begin{prop}\label{prop:convergence_lossy}
Consider Algorithm \ref{alg:robust-smart-distributed-r-admm} under Assumption \ref{ass:lossy}. Assume the pair of parameters $(\alpha, \rho)$ be such that $0<\alpha <1$ and $\rho >0$. Then, for any initial conditions, the trajectories $k \to x_i^{(i)}(k)$ and $k \to x_j^{(i)}(k)$, $i \in \mathcal{V}$, $j\in\mathcal{N}_i$ generated by the Algorithm \ref{alg:robust-smart-distributed-r-admm} converge almost surely to the optimal solution of \eqref{eq:opt_problem}, that is,
\begin{align*}
\begin{split}
	\lim_{k \to \infty} x_i^{(i)}(k) = x_i^* \\
	\lim_{k \to \infty} x_j^{(i)}(k) = x_j^* 
\end{split}, \qquad \forall i \in \mathcal{V}, j\in\mathcal{N}_i
\end{align*}
with probability one.\oprocend
\end{prop}

\smallskip

\noindent Proving Proposition \ref{prop:convergence_lossy} is achieved by showing that the randomized partition-based ADMM of Algorithm \ref{alg:robust-smart-distributed-r-admm} conforms to the stochastic Peaceman-Rachford splitting introduced in \cite{iutzeler2013asynchronous,bianchi2016coordinate} which is provably convergent. The details are available in Appendix~\ref{app:convergence}. We highlight that allowing the packet loss probability of each edge to be in general different, the partition-based R-ADMM still conforms to the stochastic PRS framework of \cite{iutzeler2013asynchronous,bianchi2016coordinate}.

\smallskip

\begin{remark}
Note that we restrict our analysis to the case of synchronous communications and updates, with the aim of investigating the performance of the R-ADMM over faulty networks. The more realistic case of asynchronous communications will be the focus of future research.\oprocend
\end{remark}

\smallskip

\begin{remark}
Observe that both in the case of reliable communications of Proposition \ref{prop:convergence} and in the lossy scenario of Proposition \ref{prop:convergence_lossy}, the convergence is guaranteed in the same region of the parameters space $(\alpha, \rho)$. In particular Algorithm \ref{alg:smart-distributed-r-admm} and the modified version \ref{alg:robust-smart-distributed-r-admm} are shown to be provably convergent for $0<\alpha<1$ and $\rho>0$. This is however only a sufficient result and the convergence might hold also in a larger region of the parameter space. Indeed this is verified by the simulations results described in \ref{sec:simulation} for the case of quadratic cost functions. Moreover, 
despite what suggested by the intuition, the larger the packet
loss probability $p$, the larger the region of convergence (though only slightly).
However, this increased region of stability is counterbalanced
by a slower convergence rate of the algorithm.
\oprocend
\end{remark}

\smallskip

\begin{remark}
It is worth to remarking that the convergence of the partition-based ADMM can be proved, both in the perfect communication scenario and in the scenario with potential packet losses, also with a time-varying step-size. Therefore it would be possible to design step-size choice criteria that speed up the convergence of the R-ADMM, which will be the object of future works.
\end{remark}

\section{Simulations}\label{sec:simulation}
This Section describes the simulative results obtained testing the effectiveness of the proposed partition-based Algorithm \ref{alg:robust-smart-distributed-r-admm}. In particular we are interested into the performance of the algorithm in the presence of packet losses due to unreliable communications. Our analysis is restricted to the case of quadratic cost functions defined as
\begin{align*}
	f_i\left(x_i^{(i)},\{x_j^{(i)}\}_{j\in\mathcal{N}_i}\right) & = \norm{A_{ii}x_i^{(i)}+\sum_{j\in\mathcal{N}_i}A_{ij}x_j^{(i)}-b_i}_{Q_i}^2
\end{align*}
where $A_{ii}\in\mathbb{R}^{r\times n}$, $A_{ji}\in\mathbb{R}^{r\times n}$ for any $j\in\mathcal{N}_i$, $b_i\in\mathbb{R}^r$, and $Q_i\in\mathbb{R}^{r\times r}$ is a symmetric and positive definite matrix. Notice that the weighted norm used in the cost function is defined as $\norm{v}_{M}^2 = v^\top M v$ for a matrix $M$ of suitable dimensions. In general the matrices $A_{ii}$ and $\{A_{ij}\}_{j\in\mathcal{N}_i}$ are different for each node, as well as the cost matrices $Q_i$. Notice that when dealing with quadratic cost functions the solution of Equation \eqref{eq:x-update-distributed} can be found in closed form.\\
Moreover, we consider the family of random geometric graphs with $N=10$ and communication radius $r=0.1$[p.u.] in which, that is, two nodes are connected if and only if their relative distance is less that $r$.\\
All the results are obtained by averaging over a set of 100 Monte Carlo runs of the simulations.

\begin{figure}[t]
      \centering
      \includegraphics[width=\columnwidth]{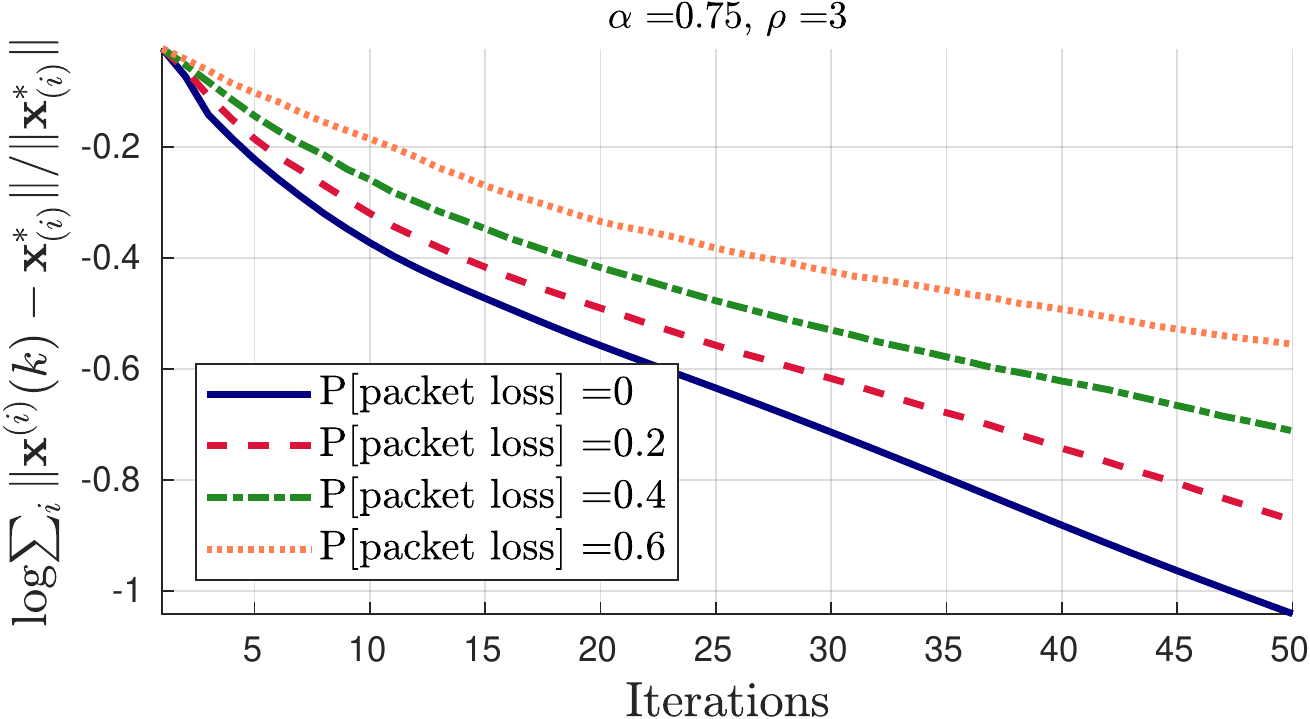}
      \caption{Evolution, in log-scale, of the relative error of Alg.~\ref{alg:robust-smart-distributed-r-admm} computed w.r.t. the unique optimal solution $\mathbf{x}^*$ as function of different values of packet loss probability $p$ for step size $\alpha=0.75$ and penalty $\rho=3$. Average over 100 Monte Carlo runs.}
      \label{fig:evolution_different_losses}
   \end{figure}

\begin{figure}[t]
      \centering
      \includegraphics[width=\columnwidth]{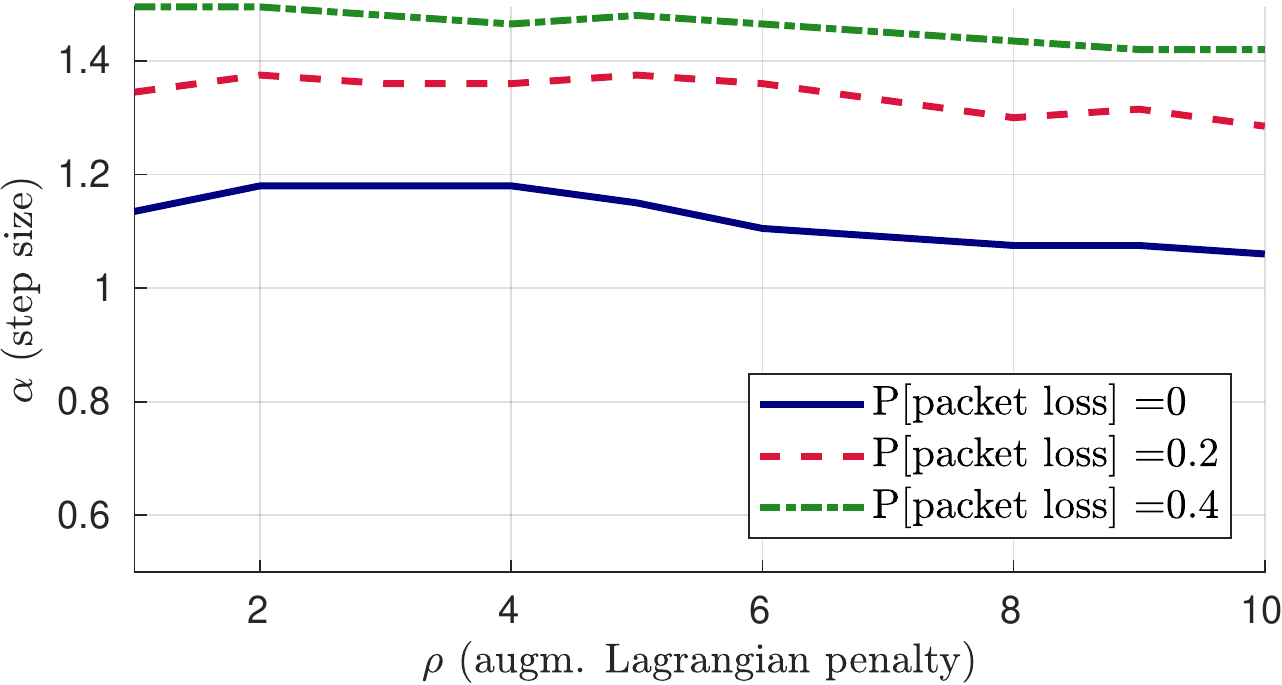}
      \caption{Stability boundaries of Alg.~\ref{alg:robust-smart-distributed-r-admm} as function of the step size $\alpha$ and the penalty $\rho$ for different values of loss probability $p$. Average over 100 Monte Carlo runs.}
      \label{fig:stability_boundaries}
   \end{figure}

In Figure~\ref{fig:evolution_different_losses} we depict the evolution of the relative error
\begin{equation}\label{eq:relative_error}
	\log \sum_{i=1}^N\frac{\norm{\mathbf{x}^{(i)}(k)-\mathbf{x}_{(i)}^*}}{\norm{\mathbf{x}_{(i)}^*}}
\end{equation}
where
\begin{equation*}
	\mathbf{x}_{(i)}^* = \begin{bmatrix}
		x_i^* \\
		\{x_j^*\}_{j\in\mathcal{N}_i}
	\end{bmatrix},
\end{equation*}
for different values of the packet loss probability, and with fixed step size $\alpha=0.75$ and penalty parameter $\rho=3$.
The presence of communication failures clearly has the effect of slowing down the convergence rate of the algorithm.\\
In Figure~\ref{fig:stability_boundaries} we report the stability boundaries of the partition-based R-ADMM for different packet loss probabilities as  functions of the tunable parameters, the step size $\alpha$ and the penalty $\rho$. In particular each curve in Figure~\ref{fig:stability_boundaries} represents the numerical boundary below which the algorithm is found to be convergent, and above which it is found to be divergent.
The results are quite interesting and will be a direction of future investigation. As was expected from the convergence result of Proposition~\ref{prop:convergence_lossy}, the convergence is guaranteed for any value of the penalty parameter $\rho$. However as the packet loss probability increases, the stability region with respect to the step size broadens. Therefore, somewhat counterintuitively, the greater the probability $p$ is, the larger the stability regions in the $(\rho,\alpha)$ space are; however this phenomenon is balanced by slower convergence rates, as depicted in Figure~\ref{fig:evolution_different_losses}.

\begin{figure}[ht]
      \centering
      \includegraphics[width=\columnwidth]{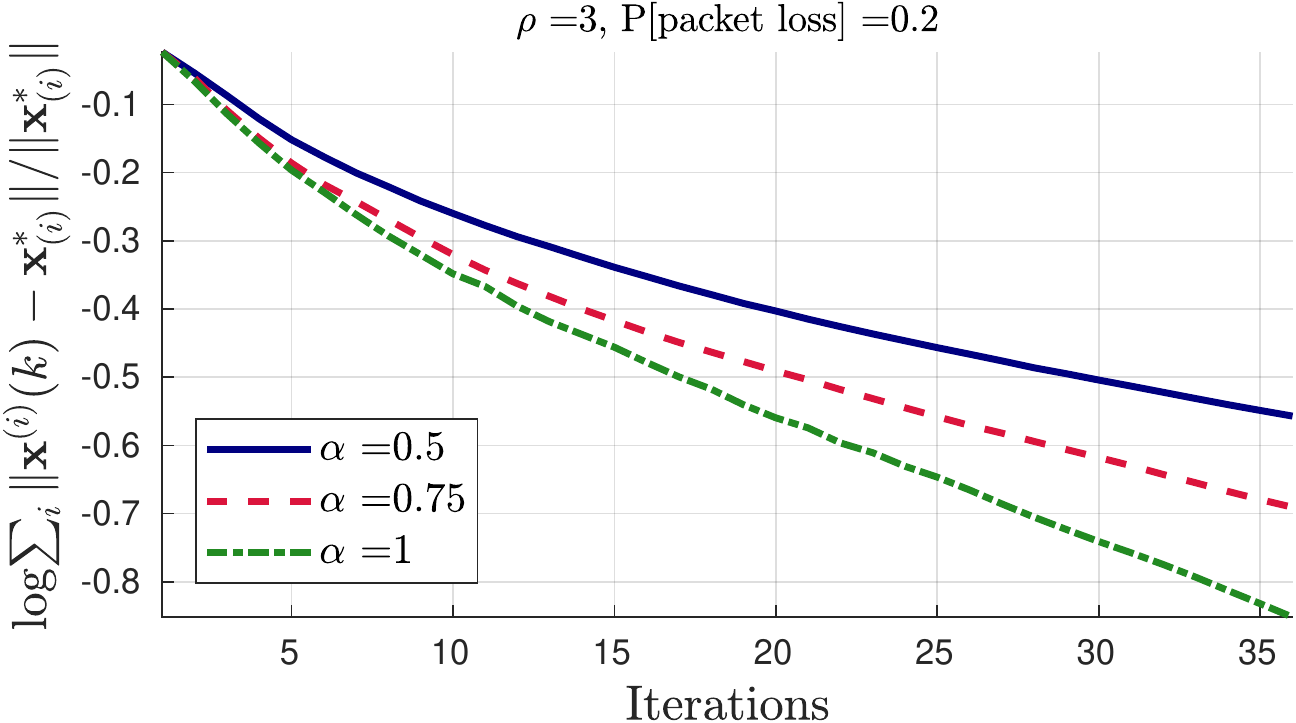}
      \caption{Evolution, in log-scale, of the relative error of Alg.~\ref{alg:robust-smart-distributed-r-admm} computed w.r.t. the unique optimal solution $\mathbf{x}^*$ as function of different values of the step size $\alpha$, with fixed packet loss probability $p=0.2$ and penalty $\rho=3$. Average over 100 Monte Carlo runs.}
      \label{fig:evolution_different_stepsizes}
   \end{figure}

\begin{figure}[ht]
      \centering
      \includegraphics[width=\columnwidth]{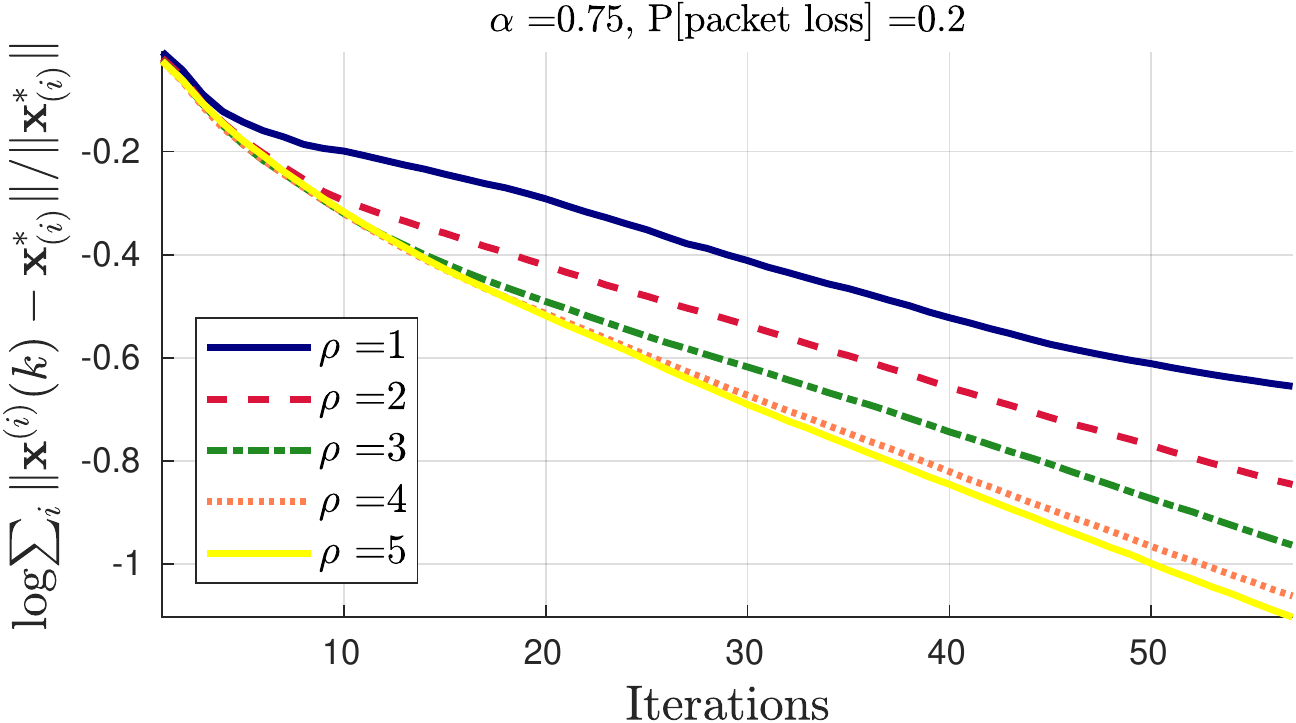}
      \caption{Evolution, in log-scale, of the relative error of Alg.~\ref{alg:robust-smart-distributed-r-admm} computed w.r.t. the unique optimal solution $\mathbf{x}^*$ as function of different values of the penalty $\rho$, with fixed packet loss probability $p=0.2$ and step size $\alpha=0.75$. Average over 100 Monte Carlo runs.}
      \label{fig:evolution_different_penalties}
   \end{figure}

The role of the tunable parameters is investigated next.\\
Figure~\ref{fig:evolution_different_stepsizes} represents the evolution of the relative error \eqref{eq:relative_error} for different values of the step size $\alpha$ and with fixed packet loss probability $p=0.2$ and penalty $\rho=3$. The use of the relaxed ADMM, instead of the classic ADMM which coincides with the R-ADMM for $\alpha=0.5$, clearly can be beneficial for the speed of convergence. In particular, inside the convergence region guaranteed by Proposition~\ref{prop:convergence_lossy}, that is $0<\alpha<1$, the rate of convergence results to be larger for values of the step size that are larger than $1/2$.
Finally, Figure~\ref{fig:evolution_different_penalties} depicts the relative error \eqref{eq:relative_error} for different values of the penalty parameter $\rho$, with step size set to $\alpha=0.75$ and packet loss probability to $p=0.2$. Recall that by Proposition~\ref{prop:convergence_lossy} the convergence of Algorithm~\ref{alg:robust-smart-distributed-r-admm} is guaranteed when the condition $\rho>0$ is satisfied. Therefore Figure~\ref{fig:evolution_different_penalties} shows that it is possible to make use of the penalty to speed up the convergence rate of the algorithm.

\section{Conclusions and Future Directions}\label{sec:conclusions}
In this paper we have presented a formulation of the relaxed ADMM tailored to distributed convex optimization with partition-based cost functions, that is, the local cost stored by each node depends on both its own state and the states of its neighbors. We have first introduced the framework describing the partition-based scenario, and then we have showed how it is possible to reformulate it so that the R-ADMM can be properly applied. Moreover we have presented an implementation of the R-ADMM that has lower memory and communications requirements than the implementation derived from a straightforward application of the algorithm.\\
The formulation of the partition-based R-ADMM that we have introduced turns out to be provably robust to random communication failures. In particular, we have rigorously proved that the region of convergence of the algorithm in the lossy scenario does not deteriorate compared to the case of reliable communications.  An interesting numerical result shows that the presence of communication failures with larger probability increases the size of the region of convergence in the space of the tunable parameters. The role of the step size and the penalty has been analyzed as well.\\
Future research will deal with to the analysis of the asynchronous case, and with the rigorous mathematical characterization of the convergence regions.

\appendices

\section{Proof of Proposition \ref{prop:distributed-r-admm}}\label{app:distributed-admm}
As we showed in Section \ref{subsec:distributed_problem} of the main paper, it is possible to reformulate the partition-based problem \eqref{eq:distributed-primal} so that it conforms to problem
\begin{align}\label{eq:primal-problem-app}
\begin{split}
	& \min_{\mathbf{x}}\ \{ f(\mathbf{x}) + \iota_{(I-P)}(\mathbf{y}) \} \\
	& \text{s.t.}\ \ A\mathbf{x}+\mathbf{y} = 0
\end{split}
\end{align}
to which the R-ADMM can be applied. The three update equations \eqref{eq:r-admm-1}, \eqref{eq:r-admm-2} and \eqref{eq:r-admm-3} that characterize the R-ADMM applied to problem \eqref{eq:primal-problem-app} yield
\begin{align}
\begin{split}
	\mathbf{y}(k+1)&=\argmin_\mathbf{y}\{\mathcal{L}_\rho(\mathbf{x}(k),\mathbf{y};\mathbf{w}(k))\\&+\rho(2\alpha-1)\langle \mathbf{y},(A\mathbf{x}(k)+\mathbf{y}(k))\rangle\}\label{eq:r-admm-1-appendix}
\end{split}\\
\begin{split}
	\mathbf{w}(k+1)&=\mathbf{w}(k)-\rho(A\mathbf{x}(k)+\mathbf{y}(k+1))\\&-\rho(2\alpha-1)(A\mathbf{x}(k)+\mathbf{y}(k))\label{eq:r-admm-2-appendix}
\end{split}\\
	\mathbf{x}(k+1)&=\argmin_\mathbf{x}\mathcal{L}_\rho(\mathbf{x},\mathbf{y}(k+1);\mathbf{w}(k+1))\label{eq:r-admm-3-appendix}
\end{align}
where $\mathbf{w}$ is the vector of Lagrange multipliers and the augmented Lagrangian is
\begin{align*}
	\mathcal{L}_\rho(\mathbf{x},\mathbf{y};\mathbf{w}) &= f(\mathbf{x}) + \iota_{(I-P)}(\mathbf{y}) -\mathbf{w}^\top(A\mathbf{x}+\mathbf{y}) \\
	& +\frac{\rho}{2}\norm{A\mathbf{x}+\mathbf{y}}^2.
\end{align*}
However, as shown in \cite{davis2016convergence}, the R-ADMM for problem \eqref{eq:primal-problem-app} can be equivalently characterized with the set of four iterates
\begin{align}
	\mathbf{y}(k)&=\argmin_{\mathbf{y}=P\mathbf{y}} \left\{-\mathbf{z}^\top(k)\mathbf{y}+\frac{\rho}{2}\|\mathbf{y}\|^2\right\} \label{eq:y-update} \\
	\mathbf{w}(k)&= \mathbf{z}(k)-\rho\mathbf{y}(k) \label{eq:w-update} \\
	\mathbf{x}(k)&=\argmin_\mathbf{x}\Big\{f(\mathbf{x})-(2\mathbf{w}(k)-\mathbf{z}(k))^\top A\mathbf{x} \nonumber \\
	&\qquad\qquad\qquad +\frac{\rho}{2}\|A\mathbf{x}\|^2\Big\} \label{eq:x-update} \\
	\mathbf{z}(k+1)&=(1-2\alpha)\mathbf{z}(k)+2\alpha(\mathbf{w}(k)-\rho A\mathbf{x}(k)) \label{eq:z-update}.
\end{align}
Similarly to what has been done in \cite{bastianello2018distributed}, it is now possible to leverage the distributed nature of problem \eqref{eq:primal-problem-app} in order to simplify Equations \eqref{eq:y-update}--\eqref{eq:z-update}.

First of all, solving the system of KKT conditions for \eqref{eq:y-update} yields $\mathbf{y}(k)=(I+P)\mathbf{z}(k)/(2\rho)$, and therefore Equations \eqref{eq:y-update}--\eqref{eq:z-update} become
\begin{align}
	\mathbf{y}(k)&= (I+P)\mathbf{z}(k)/(2\rho) \\
	\mathbf{w}(k)&= (I-P)\mathbf{z}(k)/2 \\
	\mathbf{x}(k)&=\argmin_\mathbf{x}\Big\{f(\mathbf{x})+(P\mathbf{z}(k))^\top A\mathbf{x} +\frac{\rho}{2}\|A\mathbf{x}\|^2\Big\} \label{eq:x-update-bis} \\
	\mathbf{z}(k+1)&=(1-\alpha)\mathbf{z}(k)-\alpha P\mathbf{z}(k)-2\alpha\rho A\mathbf{x}(k) \label{eq:z-update-bis}.
\end{align}
Since we are interested in the trajectory $k\to\mathbf{x}(k)$ and by the fact that the update \eqref{eq:x-update-bis} depends only on the vector $\mathbf{z}(k)$, then the R-ADMM for problem \eqref{eq:primal-problem-app} can be described by Equations \eqref{eq:x-update-bis} and \eqref{eq:z-update-bis} only.

Notice now that the trajectory $k\to\mathbf{x}(k)$ generated by \eqref{eq:x-update-bis} is equivalent to that generated by \eqref{eq:r-admm-3-appendix} if the initial condition for $\mathbf{x}$ is the same and if $\mathbf{z}(0)=\mathbf{w}(0)+\rho\mathbf{y}(0)$ since Equation \eqref{eq:w-update} has to hold at time $k=0$. Therefore Propositon \ref{prop:distributed-r-admm} is proved if we can show that \eqref{eq:x-update-bis} and \eqref{eq:z-update-bis} can be rewritten as \eqref{eq:x-update-distributed} and \eqref{eq:z-update-distributed}.

Recall that the permutation matrix $P$ swaps the element $z_i^{(i,j)}$ with the element $z_i^{(j,i)}$ of vector $\mathbf{z}$, and that the row of $A\mathbf{x}$ relative to the auxiliary variable $z_i^{(j,i)}$ is $-x_i^{(i)}$. Therefore it follows that
\begin{align*}
	(P\mathbf{z})^\top A\mathbf{x} &= \begin{bmatrix}
		\cdots & z_i^{(j,i)\top} & \cdots & z_i^{(i,j)\top} & \cdots
	\end{bmatrix}
	\begin{bmatrix}
		\vdots \\ -x_i^{(i)} \\ \vdots \\ -x_i^{(j)} \\ \vdots
	\end{bmatrix} \\
	&= - \sum_{i=1}^N \left\{ \sum_{j\in\mathcal{N}_i}z_i^{(i,j)\top} x_i^{(i)} + \sum_{j\in\mathcal{N}_i} z_j^{(i,j)\top}x_j^{(i)} \right\}.
\end{align*}
Moreover, for each node $i$ $x_i^{(i)}$ appears in $|\mathcal{N}_i|$ constraints and $\{x_j^{(i)}\}_{j\in\mathcal{N}_i}$, in one constraint each. Hence we have
\begin{equation*}
	\norm{A\mathbf{x}}^2 = |\mathcal{N}_i|\norm{x_i^{(i)}}^2+\sum_{j\in\mathcal{N}_i}\norm{x_j^{(i)}}^2.
\end{equation*}
Therefore Equations \eqref{eq:x-update-distributed} and \eqref{eq:z-update-distributed} can be derived from \eqref{eq:x-update-bis} and \eqref{eq:z-update-bis} using the particular structure of the problem, proving Proposition \ref{prop:distributed-r-admm}. \oprocendbis

\section{Proof of Propositions \ref{prop:convergence} and \ref{prop:convergence_lossy}}\label{app:convergence}
As was mentioned above, the partition-based problem can be reformulated as \eqref{eq:primal-problem-app} which can be solved by the application of the R-ADMM. Therefore both the convergence results of Propositions \ref{prop:convergence} and \ref{prop:convergence_lossy} follow from those of Propositions $2$ and $3$ of \cite{bastianello2018distributed}.\\
Indeed the R-ADMM is guaranteed to converge in both the loss-less and lossy scenarios as long as the step-size and penalty parameters are such that $0<\alpha<1$ and $\rho>0$. Moreover, the components of the primal variables vector, which in the partition-based case are the subvectors $\mathbf{x}^{(i)}$, are guaranteed to converge to the optimum value, that is, each variable $x_i^{(i)}$ converges to the optimum $x_i^*$.\oprocendbis

\addtolength{\textheight}{-14.95cm}

\bibliographystyle{./IEEEtran} 
\bibliography{./IEEEabrv,./references}

\end{document}